\documentclass[12pt]{amsart}

\usepackage{graphicx}
\newtheorem{thm}{Theorem}[section]

\newtheorem{lem}[thm]{Lemma}

\newtheorem{case}{Case}

\theoremstyle{definition}

\numberwithin{equation}{section}
\newcommand{\ben}{\begin{enumerate}}
\newcommand{\een}{\end{enumerate}}

\newcommand{\wt}{\widetilde}

\newcommand{\intr}{\text{int}\,}

\newcommand{\Figw}[4]{
\includegraphics[width=#1]{#2}
\caption{ #3 \label{#4} } }

\begin{document}

\title[Crosscap  number two, (1,1) knots in $S^3$]
{Crosscap number two knots in $S^3$ with (1,1) decompositions}%

\author[E. Ram\'{\i}rez]{Enrique Ram\'{\i}rez-Losada}
\address{Centro de Investigaci\'on en Matem\'aticas, A.C.\\
Guanajuato, Gto. 36240, M\'exico}
\email{kikis@cimat.mx}%

\author[L. G. Valdez]{Luis G. Valdez-S\'anchez}
\address{Department of Mathematical Sciences,
University of Texas at El Paso\\
El Paso, TX 79968, USA}
\email{valdez@math.utep.edu}%

%\thanks{}%


\date{\today}%

\dedicatory{Dedicated to Fico
on the occasion of his 60th birthday.}%
%\commby{}%
% ----------------------------------------------------------------
\begin{abstract}
M. Scharlemann has recently proved that any genus one tunnel number
one knot is either a satellite or 2-bridge knot, as conjectured by
H. Goda and M. Teragaito; all such knots admit a (1,1)
decomposition. In this paper we give a classification of the family
of (1,1) knots in $S^3$ with crosscap number two (i.e., bounding an
essential once-punctured Klein bottle).
%and conjecture that this list
%includes all tunnel number one, crosscap number two knots.
\end{abstract}

\maketitle

% ----------------------------------------------------------------

\section{Introduction}\label{intro}

H. Goda and M. Teragaito classified in \cite{tera3} the family of
non-simple genus one tunnel number one knots, and conjectured that
any genus one tunnel number one simple knot is a 2-bridge knot. This
conjecture was shown by H. Matsuda \cite{matsuda2} to be equivalent
to the statement that any genus one tunnel number one knot in $S^3$
admits a $(1,1)$ decomposition; it is in this form that M.
Scharlemann has recently settled it in
\cite{scharlemann4}.

In this paper we explore the family of crosscap number two tunnel
number one knots in $S^3$. Recall (cf.\ \cite{clark}) that a knot in
$S^3$ has crosscap number two if it bounds a once-punctured Klein
bottle but not a Moebius band; it was shown in
\cite{valdez6} that a knot $K$ has crosscap number two iff its
exterior contains a properly embedded {\it essential}
(incompressible and boundary incompressible, in the geometric sense)
once-punctured Klein bottle $F$, in which case $K$ is not a 2-torus
knot, and $F$ has integral boundary slope by
\cite{tera1}.

In contrast with genus one knots, a crosscap number two knot can
bound once-punctured Klein bottles with distinct boundary slopes;
however, as shown in
\cite{tera1,valdez6}, such knots are all satellite knots, with the
exception of the figure-8 knot and the Fintushel-Stern $(-2,3,7)$
pretzel knot. Here we restrict our attention to the family of
crosscap number two knots in $S^3$ which admit a $(1,1)$
decomposition; the special cases of tunnel number one satellite
knots, 2-bridge knots, and torus knots, are also discussed.

In order to state our main result we need to define a particular
family of $(1,1)$ knots in $S^3$. Let $S$ be a Heegaard torus of
$S^3$, and let $S\times I$ be a product regular neighborhood of $S$,
with $S$ corresponding to $S\times\{1/2\}$. An arc $\beta$ embedded
in $S\times I$ is called {\it monotone} if the natural projection
map $S\times I\rightarrow I$ is monotone on $\beta$. For $i=0,1$,
let $t_i$ be an embedded nontrivial circle in $S\times\{i\}$;
$t_i^*$ will denote a $(\pm 1,2)$ cable of $t_i$ relative to
$S\times\{i\}$; that is, $t_i^*$ is the boundary of a Moebius band
$B_i$ obtained by giving a half-twist to a thin annulus intersecting
$S\times\{i\}$ transversely in a core circle isotopic to $t_i$. Let
$R=\beta\times I$ be a rectangle in $S\times I$ such that $(B_0\cup
B_1)\cap R=(\partial B_0\cup
\partial B_1)\cap R=\partial\beta\times I$, and such that $\beta$
is a monotone arc in $S\times I$. Now let $K(t_0^*,t_1^*,R)$ be the
boundary of $B_0\cup R\cup B_1$ (see Fig.~\ref{aa15b}). With this
notation, the following theorem summarizes our main result.

\begin{figure}
\Figw{2.2in}{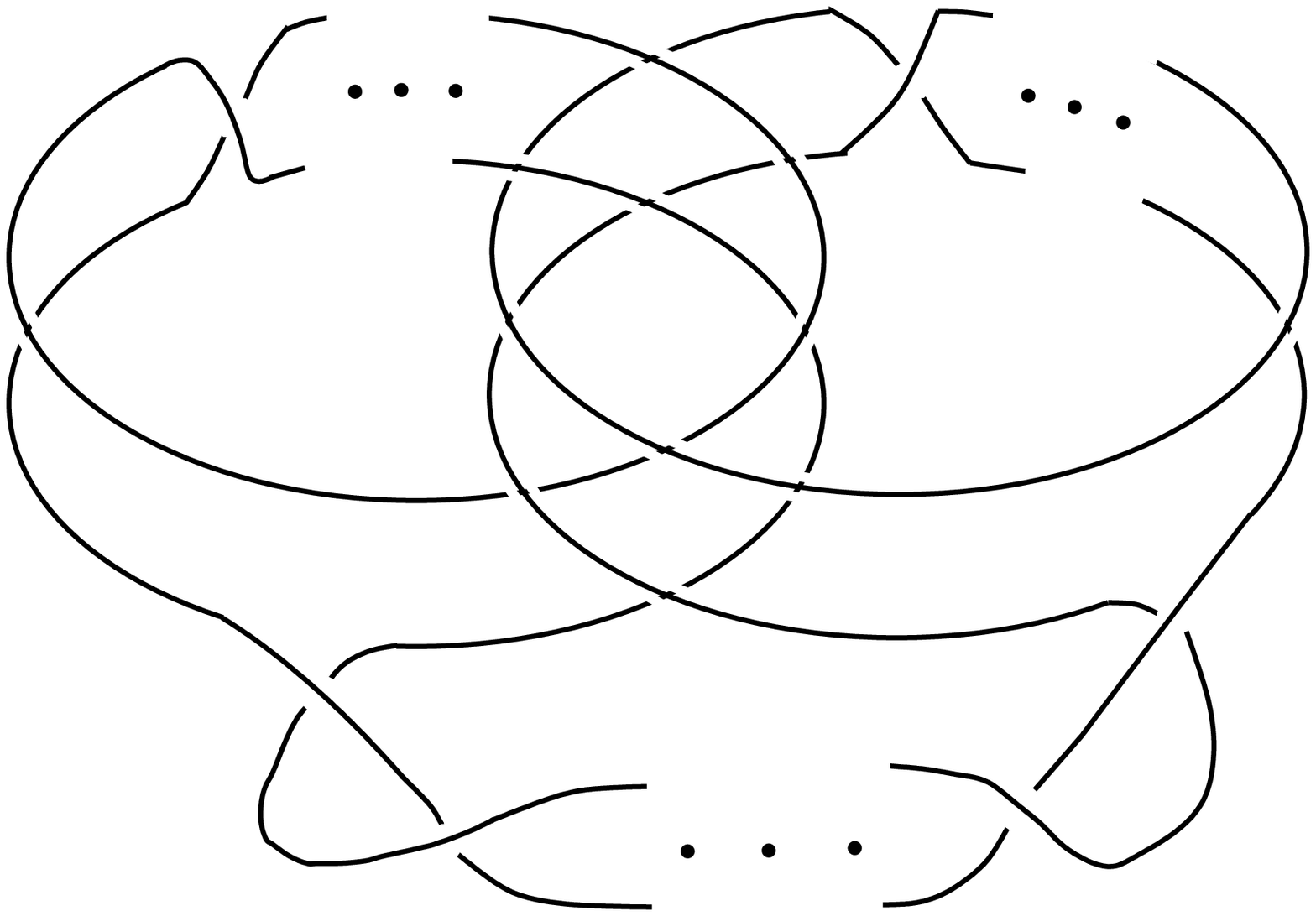}{A knot of the form
$K(t_0^*,t_1^*,R)$.}{aa15b}
\end{figure}

\begin{thm}\label{thm1}
Let $K$ be a crosscap number two knot in $S^3$. If $K$ admits a
$(1,1)$ decomposition, then $K$ is either a torus knot, a 2-bridge
knot, a satellite knot, or a knot of the form $K(t_0^*,t_1^*,R)$.
\end{thm}

The families of 2-bridge knots and tunnel number one satellite
knots, both of which admit $(1,1)$ decompositions, are of
independent interest, and we classify those having crosscap number
two explicitly; we note here (see Section~\ref{satellites}) that the
exterior $X_K$ of a tunnel number one satellite knot $K\subset S^3$
can be decomposed as the union $X_L\cup_T X_{K_0}$ for some 2-bridge
link $L$ and torus knot $K_0$. We call any  $(p,q)$ torus knot with
$|p|=2$ or $|q|=2$ a {\it $2$-torus knot}.

\begin{thm}\label{thm2}
Let $K$ be a crosscap number two knot in $S^3$; then,
\ben
\item[(a)] $K$ is a 2-bridge knot iff $K$ is a plumbing of an
annulus and a Moebius band, i.e., iff $K$ is of the form
$(2m(2n+1)-1)/(2n+1)$ for $m\neq 0$ (see Fig.~\ref{aa11}(a));\\

\item[(b)] $K$ is a tunnel number one satellite knot, with
$X_K=X_L\cup_T X_{K_0}$, iff, for some integer $m$, either\\
\ben
\item[(i)] $K_0$ is any nontrivial torus knot and
$L$ is the $4(4m+2)/(4m+1)$ or $8(m+1)/(4m+3)$ 2-bridge link
(see Fig.~\ref{aa11}(b),(c)),\\

\item[(ii)] $K_0$ is any nontrivial 2-torus knot and
$L$ is the $(8m+6)/(2m+1)$ 2-bridge link (see Fig.~\ref{aa11}(d)).
\een \een
\end{thm}

\begin{figure}
\Figw{\the\textwidth}{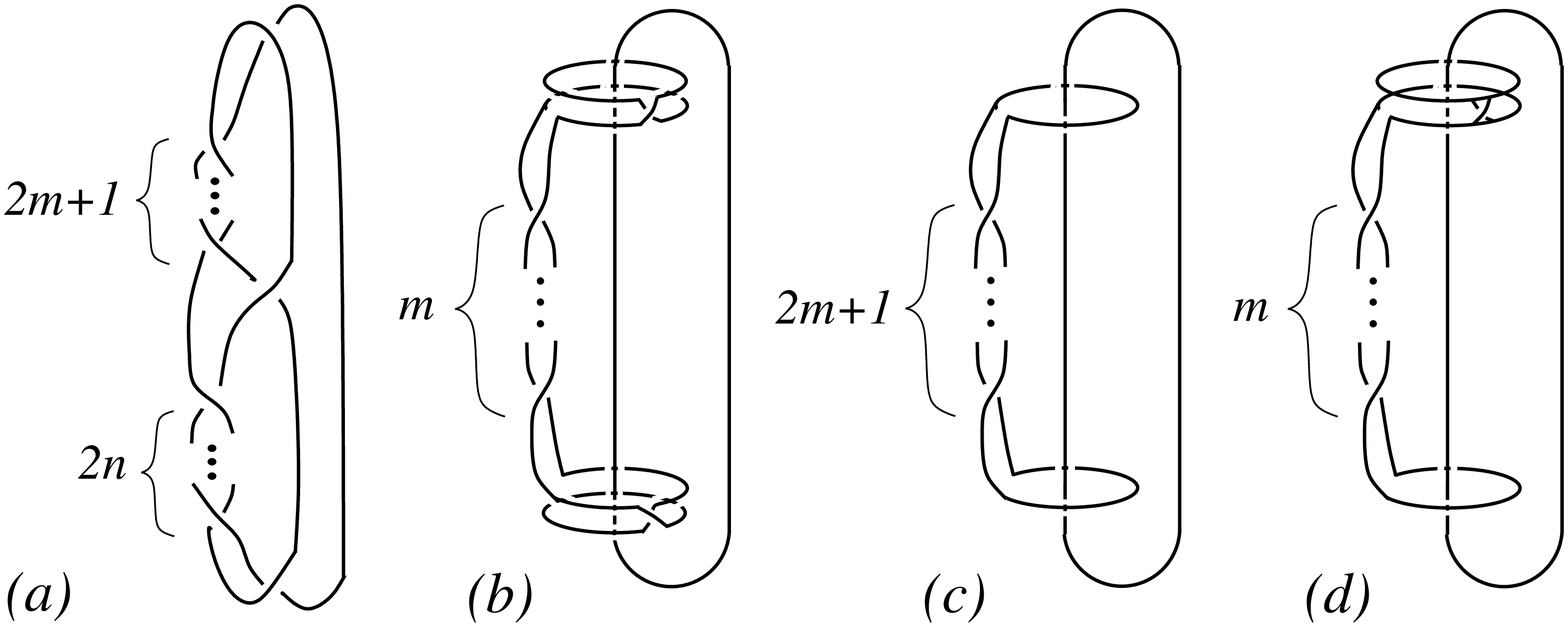}{The knots and links of
Theorem~\ref{thm2}; the integers denote half-twists.
%The $\frac{2m(2n+1)-1}{2n+1},\frac{4(4m+2)}{4m+1},
%\frac{8(m+1)}{4m+3}$, and $\frac{8m+6}{2m+1}$ 2-bridge knot and
%links.
}{aa11}
\end{figure}

We remark that in Theorem~\ref{thm2}(b) the knot $K$ is an iterated
torus knot iff $m=0$; in such case, combining the classifications of
crosscap number two cable knots in
\cite{valdez7} and of tunnel number one cable knots in
\cite{eudave3}, it follows that $K$ must be an iterated torus knot
of the form $[(4pq\pm 1,4),(p,q)]$ or $[(6p\pm 1,3),(p,2)]$ for some
integers $p,q$.

Examples of $(1,1)$ knots of the form $K(t_0^*,t_1^*,R)$ which are
neither torus, 2-bridge, nor satellites are provided by the
$(p,q,\pm 2)$ pretzel knots  with $p,q$ odd integers distinct from
$\pm 1$, as shown in Fig.~\ref{aa16}; in fact, by \cite{sakuma2},
these are the only tunnel number one pretzel knots which are not
2-bridge.

\begin{figure}
\Figw{4.3in}{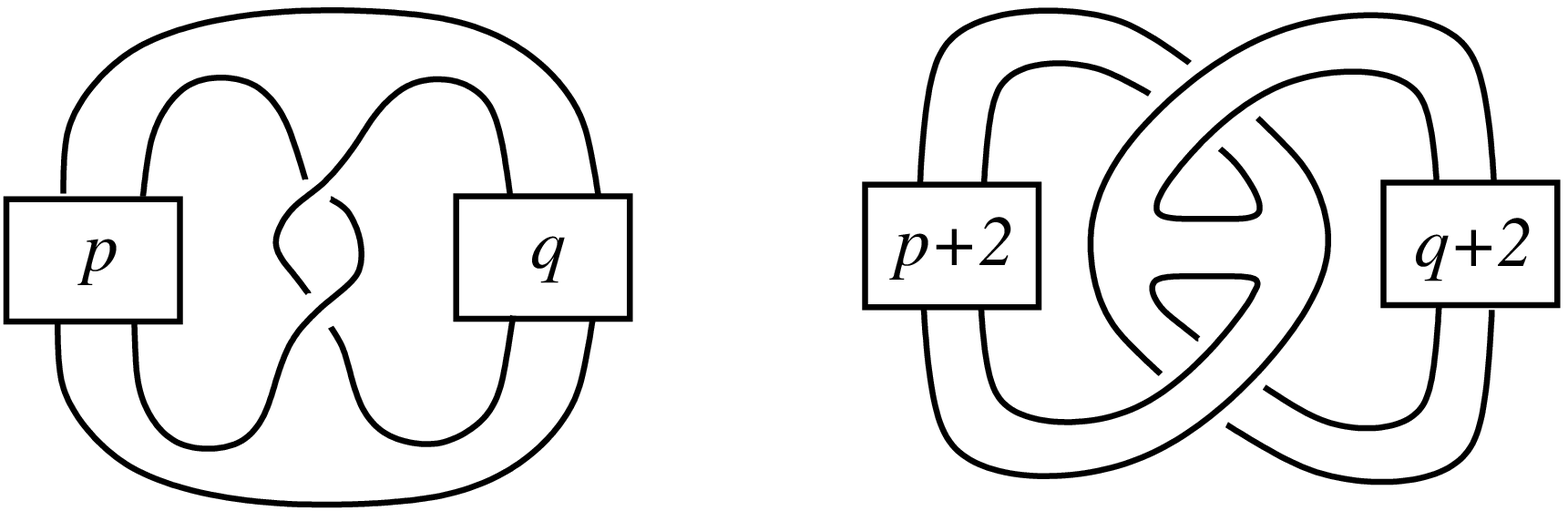}{The $(p,q,2)$ pretzel knot.}{aa16}
\end{figure}

Finally, the crosscap number two torus knots are also classified in
\cite{valdez7}; the crosscap number of torus knots in general are
determined in \cite{tera4}.

\begin{thm}[\cite{valdez7}]\label{thm3}
A $(p,q)$ torus knot has crosscap number two iff $(p,q)$ or $(q,p)$
is of the form $(3,5)$, $(3,7)$, or $(2(2m+1)n\pm 1,4n)$ for some
integers $m,n, \ n\neq 0$.\hfill\qed
\end{thm}

We will work in the smooth category. In Section~\ref{twobridge} we
discuss $(g,n)$ decompositions for knots in $S^3$ and prove
Theorem~\ref{thm2}(a). This first case, involving 2-bridge knots,
has a pleasant solution arising directly from the classification of
$\pi_1$-injective surfaces in 2-bridge knot exteriors by Hatcher and
Thurston \cite{thurs4}; we will follow and extend the basic ideas of
\cite{floyd1,floyd2,thurs4} to handle the remaining cases along
similar lines, via Morse position of essential surfaces relative to
a Heegaard surface product structure. In the process it becomes
necessary to deal with essential surfaces $\Sigma$ in knot or link
exteriors, all satisfying $\chi(\Sigma)=-1$.
Section~\ref{satellites} improves slightly on the theme of
\cite{floyd2} to allow for nonorientable essential surfaces in a
2-bridge link exterior; this is the content of Lemma~\ref{eudave},
which leads to the proof of Theorem~\ref{thm2}(b). To handle the
case of knots with a $(1,1)$ decomposition along the same lines it
is necessary to prove a statement similar to Lemma~\ref{eudave};
this is done in Section~\ref{unouno}, where Lemma~\ref{eudave2} is
established and which, along with some results from \cite{matsuda2},
leads to a proof of Theorem~\ref{thm1}. Since a once-punctured torus
$\Sigma$ also satisfies $\chi(\Sigma)=-1$, the results of this paper
can be modified to obtain the classification of genus one knots in
$S^3$ with a $(1,1)$ decomposition as well.

We want to thank Mario Eudave-Mu\~noz for making his preprint
\cite{eudave4} accessible to us, which motivated the line of
argument used in Lemma~\ref{eudave2}.

\section{$(g,n)$ decompositions and 2-bridge knots}\label{twobridge}

A knot or link $L$ in $S^3$ is said to be {\it of type $(g,n)$} if
there is a genus $g$ Heegaard splitting surface $S$ in $S^3$
bounding handlebodies $H_0, H_1$ such that, for $i=0,1$, $L$
intersects $H_i$ transversely in a trivial $n$-string arc system.
Let $S\times I$ be a product regular neighborhood of $S$ in $S^3$
and let $h:S\times I\rightarrow I$ be the natural projection map. We
denote the level surfaces $h^{-1}(r)=S\times\{r\}$ by $S_r$ for each
$0\leq r\leq 1$, and assume that $S_0\subset H_0, S_1\subset H_1$,
and that $h|S\times I\cap L$ has no critical points (so $S\times
I\cap L$ consists of monotone arcs).

Let $F$ be an essential surface properly embedded in the exterior
$X_L=S^3\setminus\intr N(L)$ of $L$; such a surface can always be
isotoped in $X_L$ so that: \ben
\item[(M1)] $F$ intersects $S_0\cup S_1$ transversely; we denote
the surfaces $F\cap H_0, F\cap H_1,F\cap S\times I$ by $F_0,F_1,\wt
F$, respectively;

\item[(M2)] each component of $\partial F$ is either a level meridian
circle of $\partial X_L$ lying in some level set $S_r$ or it is
transverse to all the level meridians circles of $\partial X_L$ in
$S\times I$;

\item[(M3)] for $i=0,1$, any component of $F_i$ containing parts
of $L$ is a {\it cancelling disk} for some arc in $L\cap H_i$ (see
Fig.~\ref{aa01}); in particular, such cancelling disks are disjoint
from any arc of $L\cap H_i$ other than the one they cancel;

\item[(M4)] $h|\wt F$ is a Morse function with a finite set
$Y(F)$ of critical points in the interior of $\wt F$, located at
different levels; in particular, $\wt F$ intersects each noncritical
level surface transversely. \een

\begin{figure}
\Figw{3.7in}{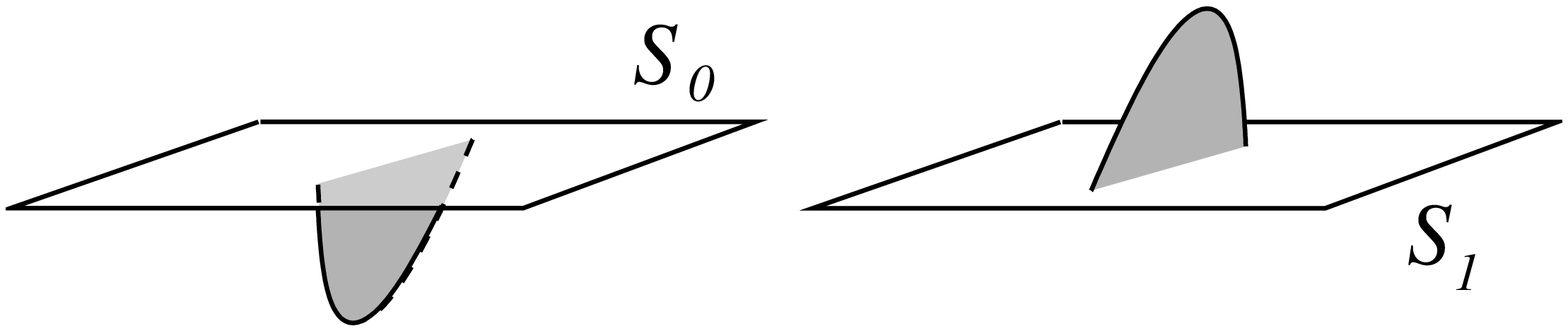}{}{aa01}
\end{figure}

We define the complexity of any surface $F$ satisfying (M1)--(M4) as
the number
$$c(F)=|\partial F_0|+|\partial F_1|+|Y(F)|,$$
where $|Z|$ stands for the number of elements in the finite set $Z$,
or the number of components of the topological space $Z$.

We say that $F$ is {\it meridionally incompressible} if whenever $F$
compresses in $S^3$ via a disk $D$ with $\partial D=D\cap F$ such
that $D$ intersects $L$ transversely in one point interior to $D$,
then $\partial D$ is parallel in $F$ to some boundary component of
$F$ which is a meridian circle in $\partial X_L$; otherwise, $F$ is
{\it meridionally compressible}. Observe that if $F$ is essential
and meridionally compressible then a `meridional surgery' on $F$
produces a new essential surface in $X_L$.

In the sequel we will concentrate in the case of knots and
2-component links $L$ of types $(0,2)$ or $(1,1)$ and certain
essential surfaces $F$ in $X_L$ with $\chi(F)=-1$. We close this
section with a proof of the first part of Theorem~\ref{thm2}.

\begin{proof}[Proof of Theorem~\ref{thm2}(a)]
Suppose $K$ is a 2-bridge knot with a $(0,2)$ decomposition relative
to some 2-sphere $S$ in $S^3$. In this context, it is proved in
\cite[Lemma 2]{thurs4} that once $F$ has been isotoped so as to
satisfy (M1)--(M4) with minimal complexity, then $F$ lies in
$S\times I$ except for the cancelling disk components of $F_0\cup
F_1$, $h|\wt F$ has only saddle critical points, $F\cap S_r$ has no
circle components for any $r$, and each saddle joins distinct level
arc components. As $\chi(F)=-1$ and $F_0$ consists of two cancelling
disks only, $h|\wt F$ has exactly three critical points, so $F$ is a
plumbing of an annulus and a Moebius band by
\cite{thurs4}; thus that $K$ must be a 2-bridge knot of the form
$(2m(2n+1)-1)/(2n+1)$ for $m\neq 0$ follows from Fig.~\ref{aa11}(a),
and the claim follows.
\end{proof}

 %%%%%%%%%%%%%%%%%%%%%%%%%%%%%%%%%%%%%%%%%%%%%%%%

\section{Satellite knots}\label{satellites}

In this section we assume that $K$ is a tunnel number one satellite
knot in $S^3$ of crosscap number two. By
\cite{eudave3,sakuma1}, the exterior $X_K=S^3\setminus\intr N(K)$
of $K$ can be decomposed as $X_L\cup_T X_{K_0}$, where
$X_L=S^3\setminus\intr N(L)$ is the exterior of some 2-bridge link
$L\subset S^3$ other than the unlink or the Hopf link and
$X_{K_0}=S^3\setminus\intr N(K_0)$ is the exterior of some
nontrivial torus knot $K_0\subset S^3$, glued along a common torus
boundary component $T$ in such a way that a meridian circle of $L$
in $T$ becomes a regular fiber of the Seifert fibration of
$X_{K_0}$.

If $F$ is any once-punctured Klein bottle, then any orientation
preserving nontrivial circle embedded in $F$ either cuts $F$ into a
pair of pants, splits off a Moebius band from $F$, or is parallel to
$\partial F$; in the first case we call such circle a {\it meridian}
of $F$, while in the second case we call it a {\it longitude} (cf.\
\cite[\S 2]{valdez6}). Notice any meridian and longitude circles of
$F$ intersect nontrivially.

As mentioned in the Introduction, $K$ has crosscap number two iff
its exterior $X_K$ contains a properly embedded essential
once-punctured Klein bottle $F$, in which case $K$ is not a 2-torus
knot and $F$ has integral boundary slope. We first show the
existence of some once-punctured Klein bottle in $X_K$ which
intersects the torus $T$ transversely in a simple way.

\begin{lem}\label{simple}
Let $K$ be a tunnel number one satellite knot in $S^3$ of crosscap
number two. Then there is an essential once-punctured Klein bottle
$F\subset X_K=X_L\cup_T X_{K_0}$ which intersects $T$ transversely
and such that either:
\ben
\item[(i)] $F$ lies in  $X_L$, or\\

\item[(ii)] $F\cap X_L$ is a once-punctured Moebius band
$F_L$ and $F\cap X_{K_0}$ is a Moebius band; in particular, $K_0$ is
a 2-torus knot. \een
\end{lem}

\begin{proof}
Let $F$ be an essential once-punctured Klein bottle in $X_K$,
necessarily having integral boundary slope; we may assume $F$ has
been isotoped so as to intersect $T$ transversely and minimally.
Hence $T\cap F$ is a disjoint collection of circles which are
nontrivial and orientation preserving in both $T$ and $F$, so each
such circle is either a meridian or longitude of $F$, or parallel to
$\partial F$ in $F$. Thus, the closure of any component of
$F\setminus T$ is either an annulus, a Moebius band, a
once-punctured Moebius band, a pair of pants, or a once punctured
Klein bottle.

Suppose $\gamma\subset T\cap F$ is a component parallel to $\partial
F$ in $F$; let $\rho$ denote the slope of a fiber of $X_{K_0}$ in
$T$. Then the component of $F\cap X_L$ containing $\partial F$ is an
annulus with the same boundary slope as $\gamma$ on $T$. If the
slope of $\gamma$ on $T$ is integral then $K$ is isotopic to $K_0$,
which is not the case; thus $\gamma$ has nonintegral slope on $T$
and so $K$ is a tunnel number one iterated torus knot with
$\Delta(\gamma,\rho)=1$ by \cite[Lemma 4.6]{eudave3}.

In particular, as any component of $F\cap X_{K_0}$ must be
incompressible and not boundary parallel in $X_{K_0}$ by minimality
of $T\cap F$, no such component can be an annulus, a Moebius band,
or a pair of pants. Therefore, $F'=F\cap X_{K_0}$ is a
once-punctured Klein bottle in $X_{K_0}$ with nonintegral boundary
slope $\gamma$ on $T$, and so, by
\cite[Lemma 4.5]{valdez6}, $F'$ must boundary compress in
$X_{K_0}$ into a Moebius band $B$ such that $\Delta(\partial
F',\partial B)=2$. But then $K_0$ is a 2-torus knot and $\partial B$
is a fiber of $X_{K_0}$, so $\Delta(\gamma,\rho)=2$, which is not
the case. Therefore, no component of $F\cap T$ is parallel to
$\partial F$ in $F$, so either $T\cap F$ is empty and (i) holds or
its components are either all meridians or all longitudes of $F$. We
now deal with the last two options.

\setcounter{case}{0}
\begin{case}
The circles $T\cap F$ are all meridians of $F$.
\end{case}

Then the component $P$ of $F\cap X_L$ containing $\partial F$ is a
pair of pants with two boundary components $c_1,c_2$ on $T$. If $A$
is an annulus in $T$ cobounded by $c_1,c_2$, then $P\cup A$ is
necessarily a once-punctured Klein bottle for $K$ which, after
pushing slightly into $X_L$, satisfies (i).

\begin{case}
The circles $T\cap F$ are all longitudes of $F$.
\end{case}

If the circles $T\cap F$ are not all parallel in $F$ then there are
two components of  $F\setminus T$  whose closures are disjoint
Moebius bands $B_1,B_2$ with boundaries on $T$. But then, if $A$ is
an annulus in $T$ cobounded by $\partial B_1,\partial B_2$, the
surface $B_1\cup B_2\cup A$ is a closed Klein bottle in $X_K\subset
S^3$, which is not possible. Hence the circles $T\cap F\subset F$
are mutually parallel in $F$, and so the component of $F\cap X_L$
which contains $\partial F$ is a once-punctured Moebius band $F_L$.
Moreover, there is a component of $F\setminus T$ whose closure is a
Moebius band $B$, properly embedded in $X_L$ or $X_{K_0}$. If $B$
lies in $X_L$ then $F\cap X_{K_0}$ is a nonempty collection of
disjoint essential annuli in $X_{K_0}$, hence $\partial B$ is the
meridian circle of a component of $L$, which implies that $B$ closes
into a projective plane in $S^3$, an impossibility. Therefore $B$
lies in $X_{K_0}$, and if $A$ is an annulus in $T$ cobounded by
$\partial F_L$ and $\partial B$ then $F_L\cup A\cup B$ can be
isotoped into a once punctured Klein bottle for $K$ satisfying (ii).
\end{proof}

Denote the components of $L$ by $K_1,K_2$, with $\partial F$
isotopic to $K_1$. We assume that a fixed 2-bridge presentation $L$
is given relative to some 2-sphere $S$ in $S^3$, and that $F$ has
been isotoped so as to satisfy (M1)--(M4) and have minimal
complexity. Notice that $H_0,H_1$ are 3-balls in this case. The next
result will be useful in the sequel.

\begin{lem}\label{eudave}
Let $\Sigma'$ be a surface in $S^3$  spanned by $K_1$ (orientable or
not) and transverse to $K_2$, such that $\Sigma=\Sigma'\cap X_L$ is
essential and meridionally incompressible in $X_L$. If $\Sigma$ is
isotoped so as to satisfy (M1)--(M4) with minimal complexity, then
$|Y(\Sigma)|=2-(\chi(\Sigma)+|\partial\Sigma|)$, and
\ben
\item[(i)]  each critical point of $h|\wt \Sigma$ is a saddle,

\item[(ii)] for $0\leq r\leq 1$ any circle component of $S_r\cap \Sigma$
is nontrivial in $S_r\setminus L$ and $\Sigma$, and

\item[(iii)] $\Sigma_0$ and $\Sigma_1$ each consists of one cancelling disk.
\een
\end{lem}

\begin{proof}
If $\Sigma$ is orientable the statement follows from the proof of
\cite[Theorem 3.1]{floyd2} without any constraints on the boundary
of $\Sigma$. If $\Sigma$ is nonorientable, the given hypothesis on
$\Sigma$ are sufficient for the arguments of \cite[Proposition
2.1]{floyd1} and \cite[Lemma 2]{thurs4} to go through and establish
(i)--(iii); the meridional incompressibility condition is needed
only for (iii), as in \cite[Theorem 3.1]{floyd2}, while the fact
that any circle component of $S_r\cap \Sigma$ is nontrivial in
$\Sigma$ follows by the argument of Lemma~\ref{eudave2}(ii). That
$|Y(F)|=2-(\chi(\Sigma)+|\partial\Sigma|)$ follows now from (i) and
(iii).
\end{proof}

\begin{proof}[Proof of Theorem~\ref{thm2}(b)]
We will split the argument into several parts, according to
Lemma~\ref{simple}.

{\bf Case (A):} \ {\it $F\subset X_L$ and $F$ is meridionally
incompressible.}

In this case Lemma~\ref{eudave} applies with $\Sigma=F$, so
$|Y(F)|=2$ and $F\cap S_0,F\cap S_1$ have no circle components. Let
$0<r_1<r_2<1$ be the levels at which the two saddles of $h|\wt{F}$
are located, and let $\alpha_0,\alpha_1$ denote the arcs $F\cap
S_0,F\cap S_1$, respectively. For any level $0<r<1$, any circle
component of $F\cap S_r$ either separates or does not separate the
points $S_r\cap K_2$; the first option is not possible by
Lemma~\ref{eudave}(ii) since $F$ is meridionally incompressible,
while in the second option it is not hard to see that, with the aid
of the cancelling disk $F_0$, $F$ compresses in $X_L$ along one such
level circle (see Fig.~\ref{aa05}).

\begin{figure}
\Figw{2in}{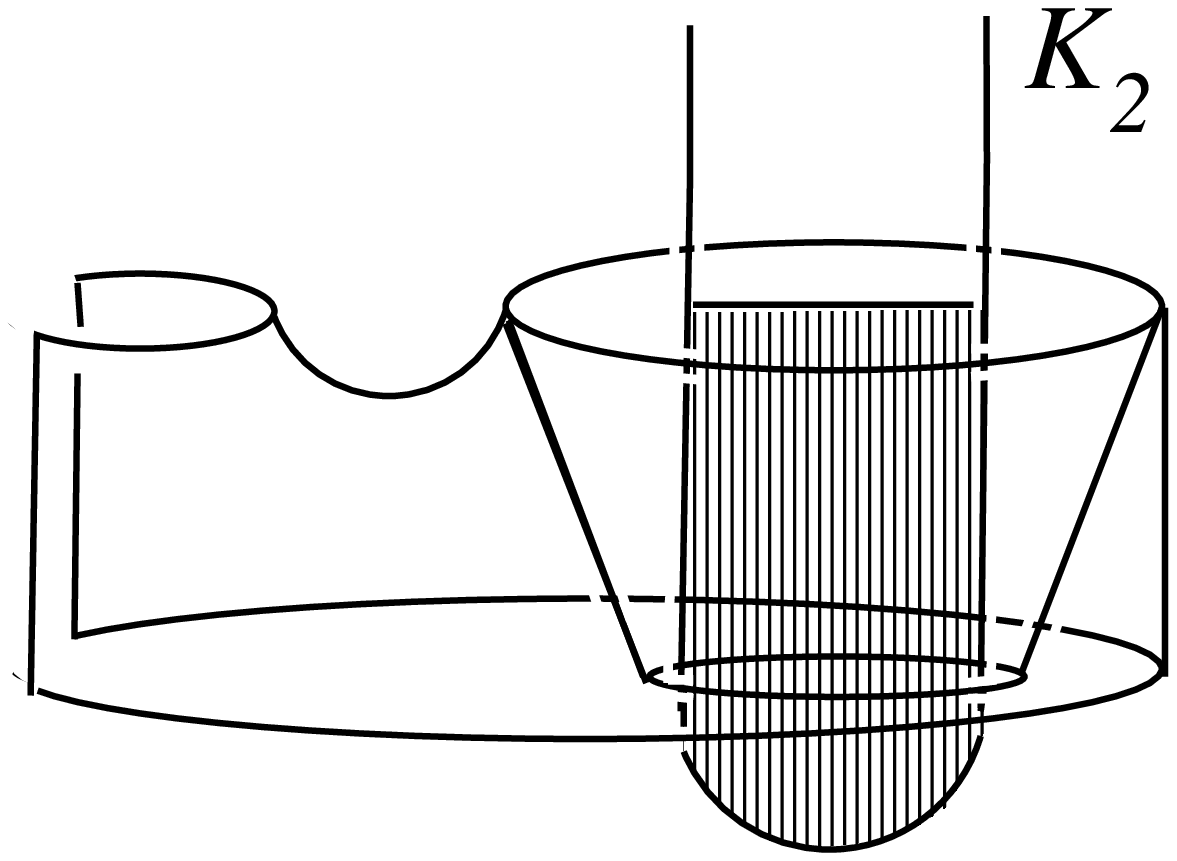}{}{aa05}
\end{figure}

Hence $S_r\cap F$ has no circle components for $0\leq r\leq 1$, so
the saddles, when seen from bottom to top and top to bottom, join
the arcs $\alpha_0,\alpha_1$, respectively, in a nonorientable
fashion (see Fig.~\ref{aa02}(a)) and so, for a sufficiently small
$\varepsilon>0$, $B_1= F\cap S\times
[r_1-\varepsilon,r_1+\varepsilon]$ and $B_2= F\cap S\times
[r_2-\varepsilon,r_2+\varepsilon]$ are Moebius bands in $F$. For
$i=1,2$, the core circle $C_i$ of $B_i$ in $S_{r_i}$ necessarily
separates the points $K_2\cap S_{r_i}$, else $C_i$ bounds a disk
$D_i$ in $S_{r_i}$ disjoint from $K_2$ as in Fig.~\ref{aa02}(b), and
a boundary compression disk for $F$ can be constructed from the
subdisk $D_i'$ of $D_i$ as in Fig.~\ref{aa02}(c); also, $\partial
B_i$ is a $(\pm 1,2)$ cable of $C_i$. Let $R$ be the rectangle $
F\cap S\times [r_1+\varepsilon,r_2-\varepsilon]\subset F$. As $h|R$
has no critical points, there exists an embedded arc $\beta$ in $R$
with one endpoint in $\partial B_1$ and the other in $\partial B_2$,
and such that $h|N(\beta)$ has no critical points for some small
regular neighborhood $N(\beta)$ of $\beta$ in $R$; thus $\beta$ is
monotone. As the once-punctured Klein bottle $ F'=B_1\cup
N(\beta)\cup B_2$ is isotopic in $X_{K_2}$ to $ F$, it follows that
the link $L$ has the form of Fig.~\ref{aa11}(a) up to isotopy (see
Fig.~\ref{aa04}), and hence that $L$ is a $4(4m+2)/(4m+1)$ 2-bridge
link. \ \qed(Case (A))

\begin{figure}
\Figw{3.5in}{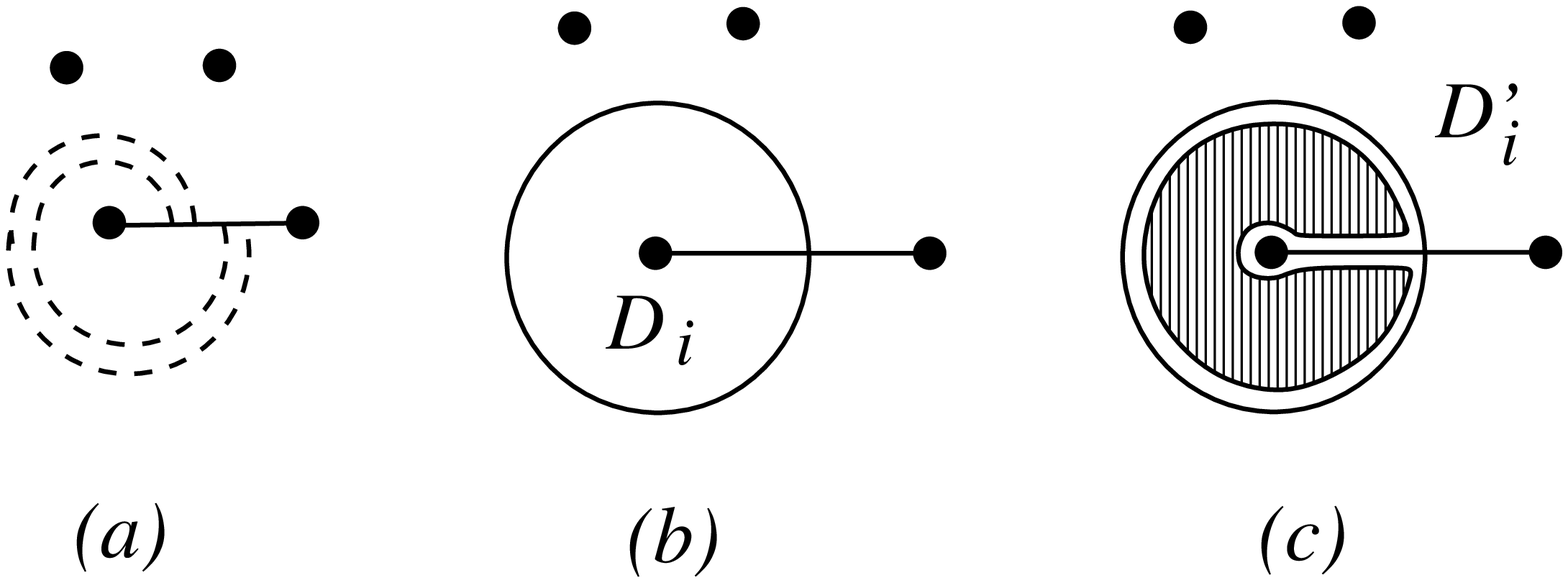}{Boundary compression of $F$.}{aa02}
\end{figure}

\begin{figure}
\Figw{3.5in}{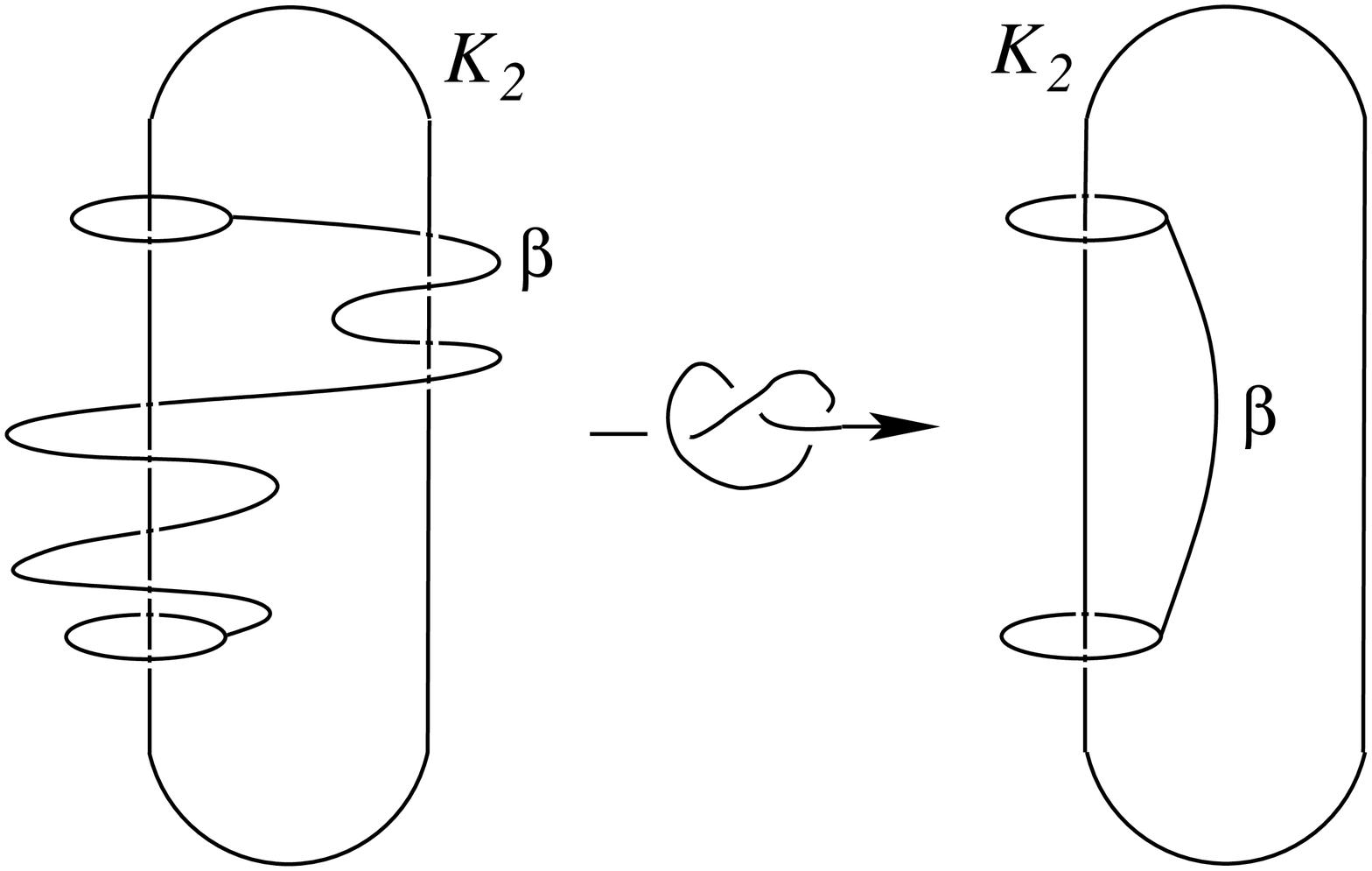}{Isotoping the arc $\beta$.}{aa04}
\end{figure}

{\bf Case (B):} \ {\it $F\subset X_L$ and $F$ is meridionally
compressible.}

Observe that if $F$  meridionally compresses along a circle
$\gamma\subset F$ then $\gamma$ must be a meridian circle of $F$:
for if $\gamma$ is trivial in $F$ then a 2-sphere in $S^3$ can be
constructed which intersects $K_2$ in one point, if $\gamma$ is a
longitude in $F$ then $S^3$ contains  $RP^2$, and if $\gamma$ is
parallel to $\partial F$ then $L$ is the Hopf link. Thus, $F$
meridionally compresses into an essential pair of pants $\Delta$ in
$X_L$, which is necessarily meridionally incompressible. By
Lemma~\ref{eudave}, we may therefore assume that $\Delta$ satisfies
(M1)--(M4) and lies within the region $S\times I$ except for the
cancelling disks $\Delta_0,\Delta_1$, and $|Y(\Delta)|=0$.

\begin{figure}
\Figw{2.7in}{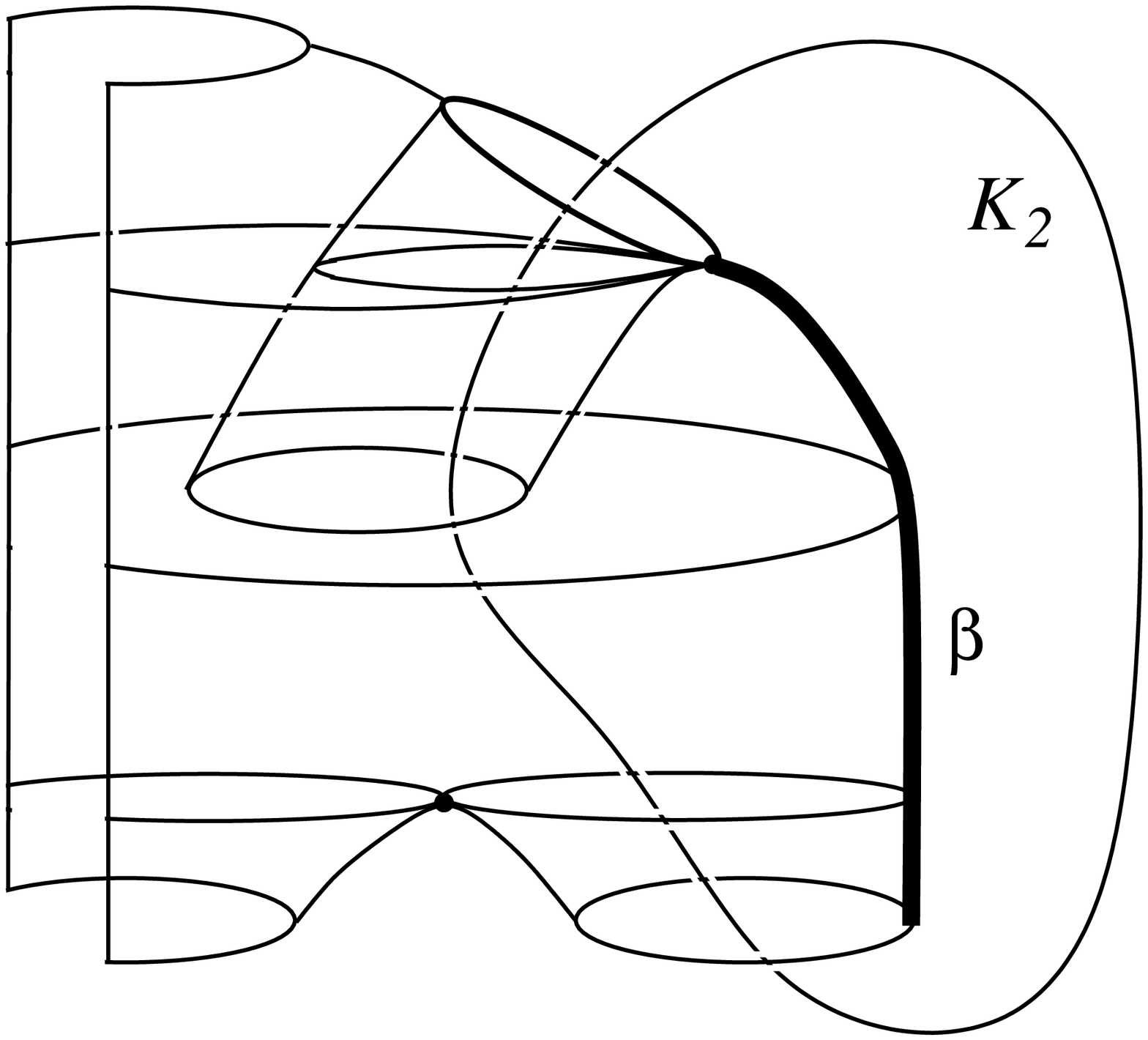}{}{aa09}
\end{figure}

Since $\Delta$ is orientable, the saddles must join the
corresponding arcs $\alpha_0=\Delta\cap S_0,\alpha_1=\Delta\cap S_1$
to themselves in an orientable fashion or to a level circle
component, when seen from bottom to top and top to bottom,
respectively. Let $C_1,C_2$ be the two level boundary circles of
$\Delta$, and let $C_3,C_4$ be the {\it limiting} circles in the
saddle levels (see Fig.~\ref{aa09}); we assume that, for $1\leq
i\leq 4$, the $C_i$'s are located at distinct levels $r_i$,
respectively. If $r_j$ and $r_k$ are the lowest and highest levels
in this list, respectively, then there exists an embedded arc
$\beta$ in $\wt\Delta$ with one endpoint in $C_j$ and the other in
$C_k$, such that $h|N(\beta)$ has no critical points for some small
regular neighborhood $N(\beta)$ of $\beta$ in $\wt\Delta$ (see
Fig.~\ref{aa09}). Then a small regular neighborhood
$N(C_j\cup\beta\cup C_k)$ in $\Delta$ yields a 2-punctured disk with
boundary isotopic to $K_1$ in $X_{K_2}$. As in Case (A), it follows
that $L$ can be isotoped into the form of Fig.~\ref{aa11}(b), so $L$
is a $8(m+1)/(4m+3)$ 2-bridge link. \
\qed(Case (B))

Therefore part (i) holds when $F\subset X_L$. We now handle the last
possible case.

{\bf Case (C):} \ {\it $F\cap X_L=F_L$.}

As for any level $0\leq r\leq 1$ each circle component of $F_L\cap
S_r$  is either parallel to the boundary circle of $F_L$ isotopic to
$K_1$, or parallel to the boundary circle of $F_L$ which is a level
meridian of $K_2$, and $L$ is neither the unlink nor the Hopf link,
it follows that $F_L$ is incompressible and meridionally
incompressible, hence Lemma~\ref{eudave} applies. Therefore, the
method of proof used in Case (B) above immediately implies that $L$
is isotopic to a link of the form of Fig.~\ref{aa11}(c), hence (ii)
holds in this case.

Since clearly any knot constructed as above has crosscap number two,
the theorem follows.
\end{proof}

\section{Knots with $(1,1)$ decompositions}\label{unouno}

In this section we assume that $K$ is a crosscap number two knot in
$S^3$ admitting a $(1,1)$ decomposition relative to some Heegaard
torus $S$ of $S^3$. In this case the handlebodies $H_0,H_1$ are
solid tori with meridian disks of slope $\mu_0,\mu_1$ in $S_0,S_1$,
respectively. For $\{i,j\}=\{0,1\}$, we project $\mu_j$ onto $S_i$,
continue to denote such projection by $\mu_j$, and frame $S_i$ via
the circles $\mu_i,\mu_j$, so that a {\it $(p,q)$-circle} in $S_i$
means a circle embedded in $S_i$ isotopic to $p\mu_i+q\mu_j$; thus
$S_i$ gets the standard framing as the boundary of the exterior of
the core of $H_i$, and a $(p,q)$-circle in $S_0$ is isotopic in
$S\times I$ to a $(q,p)$-circle in $S_1$.

Before studying the associated essential once-punctured Klein bottle
for $K$, we prove a statement similar to Lemma~\ref{eudave} in the
present context.
%Recall that $\wt{\Sigma}=\Sigma\cap S\times I$.

\begin{lem}\label{eudave2}
Suppose $K$ is not a torus knot. Let $\Sigma'$ be a spanning surface
for $K$ in $S^3$ (orientable or not) such that $\Sigma=\Sigma'\cap
X_K$ is essential in $X_K$. If $\Sigma$ is isotoped so as to satisfy
(M1)--(M4) with minimal complexity, then
$|Y(\Sigma)|=1-\chi(\Sigma)$, and \ben
\item[(i)]  each critical point of $h|\wt \Sigma$ is a saddle,

\item[(ii)] for $0\leq r\leq 1$ any circle component of $S_r\cap \Sigma$
is nontrivial in $S_r\setminus K$ and $\Sigma$, and not parallel in
$\Sigma$ to $\partial\Sigma$,

\item[(iii)] for $i=0,1$ $\Sigma_i$ consists of one cancelling
disk and either one Moebius band and some annuli components, or a
collection of disjoint annuli each having boundary slope $(p_i,q_i)$
in $S_i$ with $|q_i|\geq 2$, and

\item[(iv)] the saddle closest to either the 0-level or 1-level
does not join circle components.

\een
\end{lem}

\begin{proof}
Part (i) follows from the argument of \cite[Proposition
2.1]{floyd1}.

Suppose now that $\gamma$ is a circle component of $S_r\cap\Sigma$
for some level $0\leq r\leq 1$. If $\gamma$ bounds a disk $D$ in
$S_r\setminus K$ then $\gamma$ bounds a disk $D'$ in $\Sigma$, since
$\Sigma$ is incompressible in $X_K$. Construct a surface $\Sigma''$
isotopic to $\Sigma$ from $(\Sigma\setminus D')\cup D$ by pushing
$D$ slightly above or below $S_r$ so that $\Sigma''$ satisfies
(M1)--(M4) and the singularities of $h|\wt\Sigma''$ are exactly
those of $h|\wt\Sigma\setminus D'$ with an additional local extremum
in the interior of $D$; thus, $h|\wt\Sigma''$ has at most
$|Y(\Sigma)|+1$ critical points.

If $D'$ is disjoint from $S_0\cup S_1$ then $D'$ lies in $S\times I$
and, since $\partial D'$ is level, $h|D'$ has a local extremum in
$\intr D'$, contradicting (i). If $D'$ intersects $S_0\cup S_1$ then
$|\partial\Sigma''_0|+|\partial\Sigma''_1|<|\partial\Sigma_0|+|\partial\Sigma_1|$
while $|Y(\Sigma'')|\leq|Y(\Sigma)|+1 $, hence $c(\Sigma'')\leq
c(\Sigma)$ and so $c(\Sigma'')=c(\Sigma)$ by minimality of
$c(\Sigma)$, again contradicting (i). Therefore, $\gamma$ is
nontrivial in $S_r\setminus K$ and, since $K$ is not a torus knot,
$\gamma$ is not parallel in $\Sigma$ to $\partial\Sigma$. Thus it
only remains to verify that $\gamma$ is nontrivial in $\Sigma$ for
(ii) to hold, which we will do by the end of the proof.

If some component of $\Sigma_0$, other than the cancelling disk,
compresses in $H_0$, then there is one such component $\sigma$ which
compresses in $H_0$ via a disk $D$ disjoint from all other
components of $\Sigma_0$. Since $\Sigma$ is essential in $X_K$,
$\partial D$ bounds a disk $D'$ in $\Sigma$. Let
$\Sigma''=(\Sigma\setminus D')\cup D$. Then $h|\wt\Sigma''$ has at
most $|Y(\Sigma)|$ singular points and, since $\intr D'$ necessarily
intersects $S_0\cup S_1$,
%the closure of some disk component of $D'\setminus (S_0\cup S_1)$ lies in
%$\Sigma_0\cup\Sigma_1$ by (i), hence
$|\partial\Sigma''_0|+|\partial\Sigma''_1|<|\partial\Sigma_0|+|\partial\Sigma_1|$
and so $c(\Sigma'')<c(\Sigma)$, an impossibility. Therefore, any
component of $\Sigma_0$ is incompressible in $H_0$, hence it must be
either an annulus, a Moebius band, or a disk; since $H_0$ is a solid
torus, $\Sigma_0$ may have at most one Moebius band component.

Suppose $\Sigma_0$ has an annulus component $\sigma$; then
%$\Sigma_0$ has no Moebius band components, and
$\sigma$ separates $H_0$ into two pieces, one of which contains the
cancelling disk component of $\Sigma_0$. If $\sigma$ is parallel in
$H_0$ into $S_0$ away from all other components of $\Sigma_0$ then
$\sigma$ can be pushed into the region $S\times I$; notice this is
the case if the slope of $\sigma$ in $S_0$ is of the form
$(p_0,q_0)$ with $|q_0|=1$. It is then possible to isotope $\sigma$
and $\Sigma$ appropriately, so that $h|\sigma$ has one saddle and
one local minimum and $\Sigma$ continues to satisfy (M1)--(M4);
hence $|\partial\Sigma_0|$ will decrease by two while $|Y(\Sigma)|$
will increase by two, and so $c(\Sigma)$ will remain minimal.
However, this time $h|\wt{\Sigma}$ has a local minimum critical
point in $\sigma$, contradicting (i). Therefore, since $\sigma$ is
incompressible in $H_0$, any boundary component of $\sigma$ must be
nontrivial in $S_0\setminus K$ and distinct from $\mu_0$, so it
follows that the boundary slope of $\sigma$ in $S_0$ is of the form
$(p_0,q_0)$ with $|q_0| \geq 2$.

Consider the first saddle above level $0$; if it joins a circle
component $\gamma$ of $\Sigma\cap S_0$ to itself or to another such
circle component then it is possible to lower the saddle below level
$S_0$ while satisfying (M1)--(M4), thus reducing the value of
$c(\Sigma)$, which is not possible. Hence (iv) holds, and the first
saddle above level 0 joins the arc component $\alpha_0$ of
$S_0\cap\Sigma$ to itself or to a circle component.

Suppose now that $\sigma$ is a disk component of $\Sigma_0$ other
than the cancelling disk; then $\sigma$ is either a trivial disk or
a meridian disk of $H_0$. In the first case, $\sigma$ separates
$H_0$ into a 3-ball $B^3$ and a solid torus, with the cancelling
disk of $\Sigma_0$ contained in $B^3$ by the first part of (ii); we
may further assume that $\partial\sigma$ and $\alpha_0$ are adjacent
in $S_0$. Consider the first saddle above level 0. If it joins the
arc component $\alpha_0$ of $\Sigma\cap S_r$ to itself then either a
Moebius band is created by the saddle with core a circle bounding a
disk in the saddle level, so $\Sigma$ is boundary compressible (see
Fig.~\ref{aa02}), or a trivial circle component is created in a
level slightly above the saddle level, contradicting the first part
of (ii). If the saddle joins $\partial\sigma$ to $\alpha_0$ then
pushing down the saddle slightly below level 0 isotopes $\Sigma$ so
as to still satisfy (M1)--(M4) but lowers its complexity. Since by
(iv) these are the only possibilities for the first saddle, if
$\Sigma_0$ contains any disk components other than the cancelling
disk then all such components are meridian disks of $H_0$. The
analysis of the possible scenarios for the first saddle above level
0 is similar to that of the previous cases, except for when the
saddle joins $\alpha_0$ to itself as in Fig.~\ref{aa12}(a). In such
case, if $r$ is the level of the first saddle above level 0, the
Moebius band created by the saddle has as core a circle in $S_r$
which bounds a meridian disk of the solid torus bounded by $S_r$
below the level $S_r$ (see Fig~\ref{aa12}(b)). The situation is
similar to that of Fig.~\ref{aa02}, so $\Sigma$ is boundary
compressible, which is not the case. Hence $\Sigma_0$, and similarly
$\Sigma_1$, has no such disk components and (iii) holds.

Now let $0\leq r\leq 1$ and $\gamma$ be any circle component of
$(S_0\cup S_r\cup S_1)\cap\Sigma$.
%
%Clearly, such a circle must be orientation preserving in $\Sigma$,
%hence it must be either a trivial circle, a nontrivial circle, or
%a circle parallel to $\partial\Sigma$ in $\Sigma$. We have seen
%that the latter can not be the case since $K$ is not a torus knot
%
If $\gamma$ is trivial and innermost in $\Sigma$ then it bounds a
subdisk $D$ in $\Sigma$ with interior disjoint from $S_0\cup S_r\cup
S_1$, hence $D$ lies either in $\Sigma_0,\Sigma_1$, or $S\times I$.
But, as shown above, neither $\Sigma_0$ nor $\Sigma_1$ have disk
components other than the cancelling disks, and if $D$ lies in
$S\times I$ then, as $\partial D=\gamma$ is level, $h|D$ must have a
local extremum in $\intr D$, contradicting (i). Hence $\gamma$ is
nontrivial in $\Sigma$ and so the proof of (ii) is complete. That
$|Y(\Sigma)|=1-\chi(\Sigma)$ now follows from (i) and (iii).
\end{proof}

\begin{figure}
\Figw{5in}{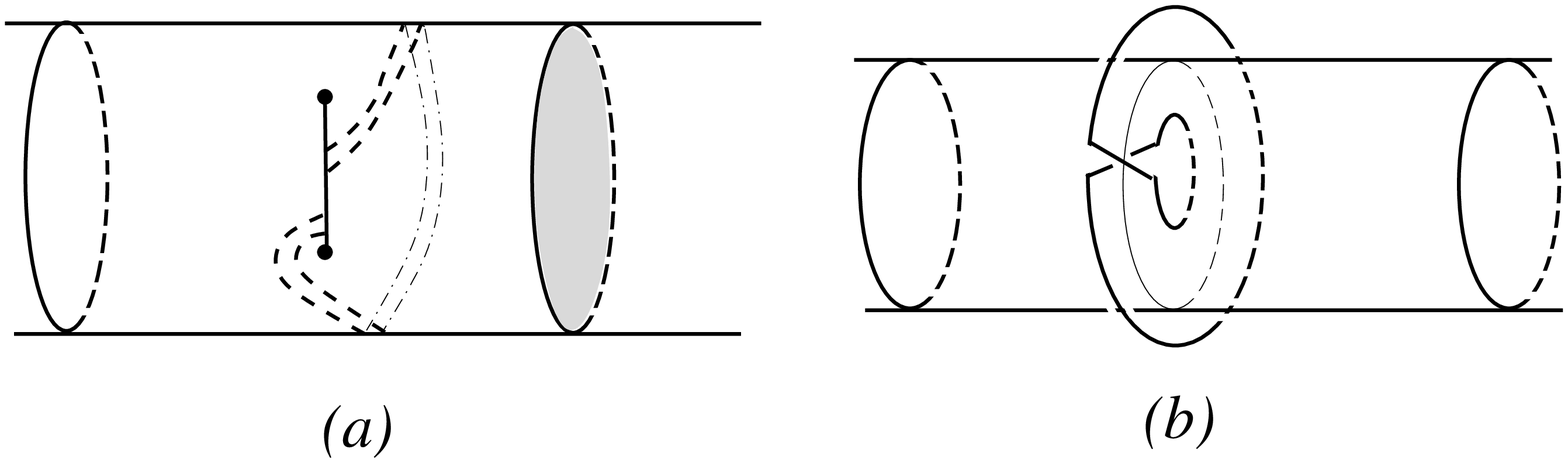}{}{aa12}
\end{figure}

In preparation for the proof of Theorem~\ref{thm1}, the following
result specializes Lemma~\ref{eudave2} to the case when $\Sigma$ is
a once punctured Klein bottle $F$; its first part is a slight
generalization of a construction by Matsuda in
\cite[pp.\ 2161--2162]{matsuda2}. We will say that an essential
annulus $A$ properly embedded in $S\times I$ is an {\it $F$-spanning
annulus} if $A$ can be isotoped so as to be disjoint from the
component of $\wt{F}=F\cap S\times I$ containing parts of $K$, and
its boundary slope in $S_0$ is of the form $(p,q)$ for some
$|p|,|q|\geq 2$. Notice that an $F$-spanning annulus $A$ is isotopic
in $S\times I$ to the annulus $(\partial A\cap S_0)\times I$, and
its boundary component in $S_1$ has slope $(q,p)$.

\begin{lem}\label{annulus}
Let $F$ be an essential once-punctured Klein bottle spanned by $K$
which has been isotoped so as to satisfy (M1)--(M4) with minimal
complexity. If there is an $F$-spanning annulus in $S\times I$
having boundary slope $(p,q)$ in $S_0$ then $K$ is either a $(p,q)$
torus knot or a satellite of a $(p,q)$-torus knot; otherwise,
$F\cap(S_0\cup S_1)$ has at most two circle components.

%besides the cancelling disk components, $F_0\cup F_1$ has at most
%two components, and if two then they must be Moebius bands.
\end{lem}

\begin{proof}
Let $F'$ denote the component of $\wt F=F\cap S\times I$ containing
parts of $K$. Let $A$ be an $F$-spanning annulus with boundary slope
$(p,q)$ in $S_0$, and suppose $K$ is not a $(p,q)$ torus knot. By
Lemma~\ref{eudave2}(ii),(iii), $F'$ is either a once-punctured
Moebius band or a pair of pants embedded in the solid torus
$V=S\times I\setminus\intr N(A)$, where $N(A)$ is a small regular
neighborhood of $A$ in $S\times I$. In either case, $\partial F'$
has one component $K'\subset \intr V$ which is isotopic to $K$ in
$S^3$, and one or two more components embedded in $\partial V$, each
running once around $V$. Notice that $V$ is a regular neighborhood
of a $(p,q)$ torus knot, so $K'$ is not a core of $V$.

If $F'$ is a once-punctured Moebius band then $K'$ is a nontrivial
knot in $V$ with odd winding number. If $F'$ is a pair of pants
then, by Lemma~\ref{eudave2}(ii),(iii), the closure of $F\setminus
F'$ consists either of two Moebius bands or an annulus with core a
meridian circle of $F$. In the first case $F_0$ and $F_1$ each have
a Moebius band component which, due to the presence of the spanning
annulus $A$, have boundary slopes $(p,q)$ and $(q,p)$, respectively,
an impossibility since then $|p|=|q|=2$; in the latter case, the
closure of the annulus $F\setminus F'$ intersects $V$ in annuli
running once around $V$, thus it can be isotoped in $S^3$, away from
$F'$, into $S^3\setminus\intr V$, and so the components of $\partial
F'$ other than $K'$ must be coherently oriented in $\partial V$;
therefore $K'$ has winding number two in $V$ and hence it is a
nontrivial satellite of the core of $V$. The first part of the lemma
follows.

Suppose now that $F\cap(S_0\cup S_1)$ has at least three circle
components; if, say, three such components lie in $S_0$, or at least
two lie in $S_0$ and at least one in $S_1$, then, since $|Y(F)|=2$
by Lemma~\ref{eudave2} and the saddles do not join circle
components, at least one of the circle components of $F\cap S_0$
must flow along an annulus component of $\wt{F}$ from $S_0$ to $S_1$
without interacting with the saddles. Thus $\wt{F}$ has at least one
annulus component which, by Lemma~\ref{eudave2}(iii), has boundary
slope of the form $(p,q)$ in $S_0$ for some $|p|,|q|\geq 2$, and so
must be an $F$-spanning annulus. Thus the second part of the lemma
follows.
\end{proof}

\begin{proof}[Proof of Theorem~\ref{thm1}]
Let $K$ be a crosscap number two knot in $S^3$, and let $F$ be an
essential once-punctured Klein bottle spanned by $K$; we assume $F$
has been isotoped so as to satisfy (M1)--(M4) with minimal
complexity. To simplify notation, let $F_0',F_1'$ denote the
components of $F_0,F_1$, respectively, other than the cancelling
disks. By Lemma~\ref{annulus}, if $K$ is neither a torus nor a
satellite knot then $S\times I$ contains no $F$-spanning annuli and
$F\cap(S_0\cup S_1)$ has at most two circle components; thus,
without loss of generality, $F'_0$ and $F'_1$ fit in one of the
following cases.

%{\bf Case (A):} \ {\it $F_0'$ is an annulus and $F_1'$ is a
%Moebius band.}
%
%By Lemma~\ref{eudave2}(iv), one of the circle components of the
%annulus $\partial F_0'$ must survive the interaction with the
%first saddle above the 0-level; on the other hand, the only circle
%component of $\partial F_1'$ disappears after interacting with the
%first saddle below the 1-level. This is impossible since $Y(F)=2$
%and the number of circle components of $F\cap S_r$ must be
%constant between consecutive saddle levels.

{\bf Case (A):} \ {\it $F_0'$ is an annulus and $F_1'$ is empty.}

Fig.~\ref{aa14}(a) shows the only possible construction (abstractly)
of the surface $F$, starting from $F_0$, via the two saddles of
$h|\wt F$. By Lemma~\ref{eudave2}(iii), the boundary slope of the
annulus $F_0'$ in $S_0$ is of the form $(p,q)$ with $|q|\geq 2$. It
is not hard to see that the boundary circle $C$ of the annulus
$\partial F_0'$ in Fig.~\ref{aa14}(a) bounds an essential annulus
$A$ in $S\times I\setminus F$, hence $|p|=1$ since $S\times I$ has
no $F$-spanning annuli, so $K$ is a 2-bridge knot by the argument of
\cite[pp.\ 2161--2162]{matsuda2}.

\begin{figure}
\Figw{\the\textwidth}{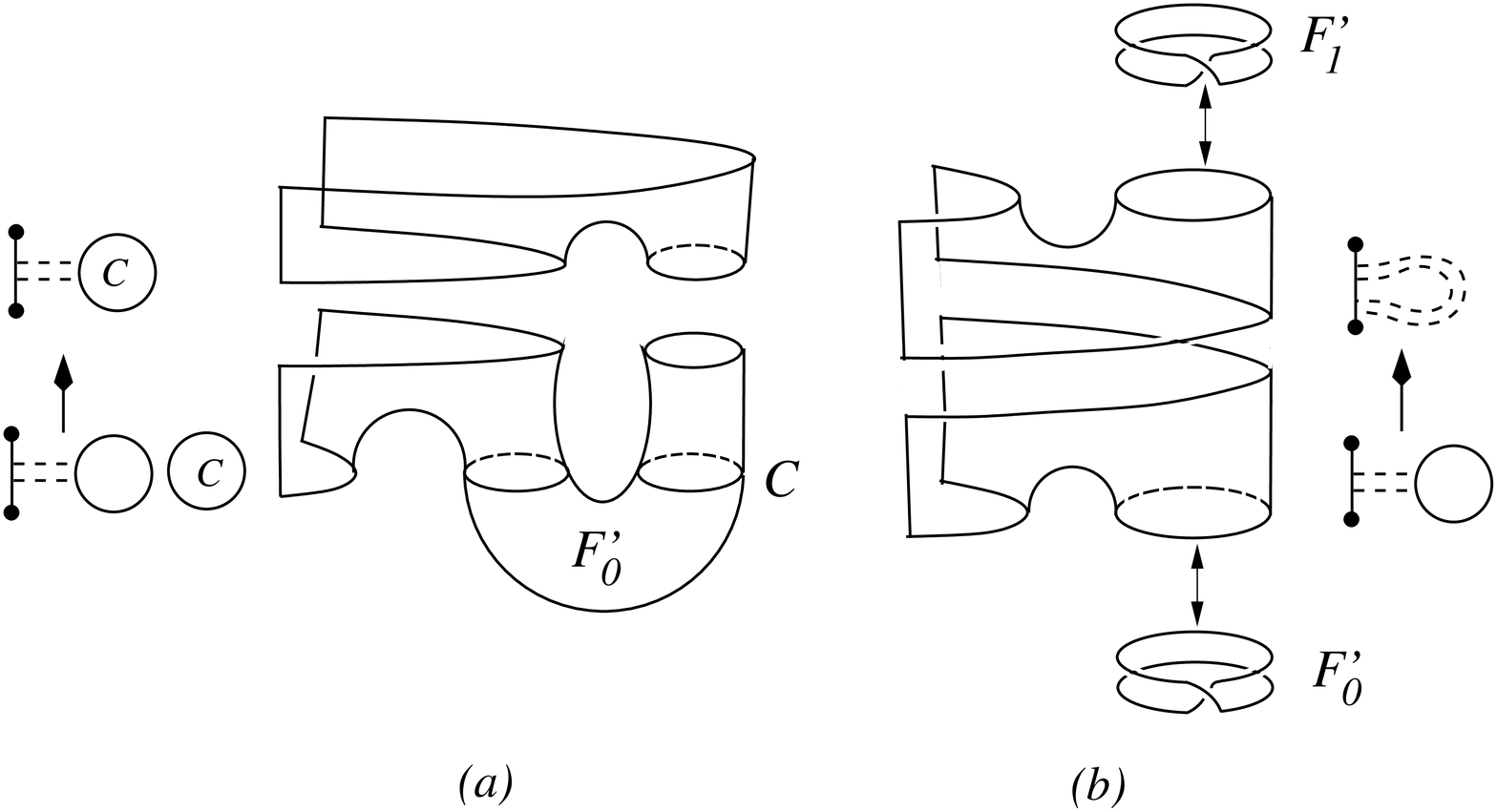}{Construction of $F$ via saddles
in Cases (A) and (B).}{aa14}
\end{figure}

{\bf Case (B):} \ {\it Both $F_0'$ and $F_1'$ are Moebius bands.}

The only possibility in this case is the one shown (abstractly) in
Fig.~\ref{aa14}(b): for otherwise, by Lemma~\ref{eudave2}(iv), the
first saddle above the 0-level would join the arc component of
$S_0\cap F$ with itself, necessarily in an orientable fashion, and
so $S_r\cap F$ would have two circle components for any level $r$ in
between the saddle levels; but then the first saddle below the
1-level must join the circle component of $S_1\cap F$ with itself,
contradicting Lemma~\ref{eudave2}(iv).

Hence $\wt F$ is a pair of pants and, since all the critical points
of $h|\wt F$ are saddles, there exists an embedded arc $\beta$ in
$\wt F$ with one endpoint in $\partial F_0'$ and the other in
$\partial F_1'$ which is monotone in $S\times I$ and such that $h|R$
has no critical points for some small regular neighborhood $R$ of
$\beta$ in $\wt F$. Observe that, for $i=0,1$, if $\partial F_i'$ is
a $(p_i,2)$-circle in $S_i$, then $F_i'$ is isotopic in $S^3$ to a
Moebius band $B_i$ which is a $(1,2)$ cable of a $(p_i,1)$-circle
$t_i$ in $S_i$. Therefore the once-punctured Klein bottle $F_0'\cup
R\cup F_1'$ can be isotoped into $B_0\cup R'\cup B_1$ for some
monotone subrectangle $R'$ of $R$. As $F_0'\cup R'\cup F_1'$ is
isotopic to $F$ in $S^3$, it follows that $K$ is a knot of the form
$K(t_0^*,t_1^*,R')$.

{\bf Case (C):} \ {\it Both $F_0'$ and $F_1'$ are empty.}

\begin{figure}
\Figw{\the\textwidth}{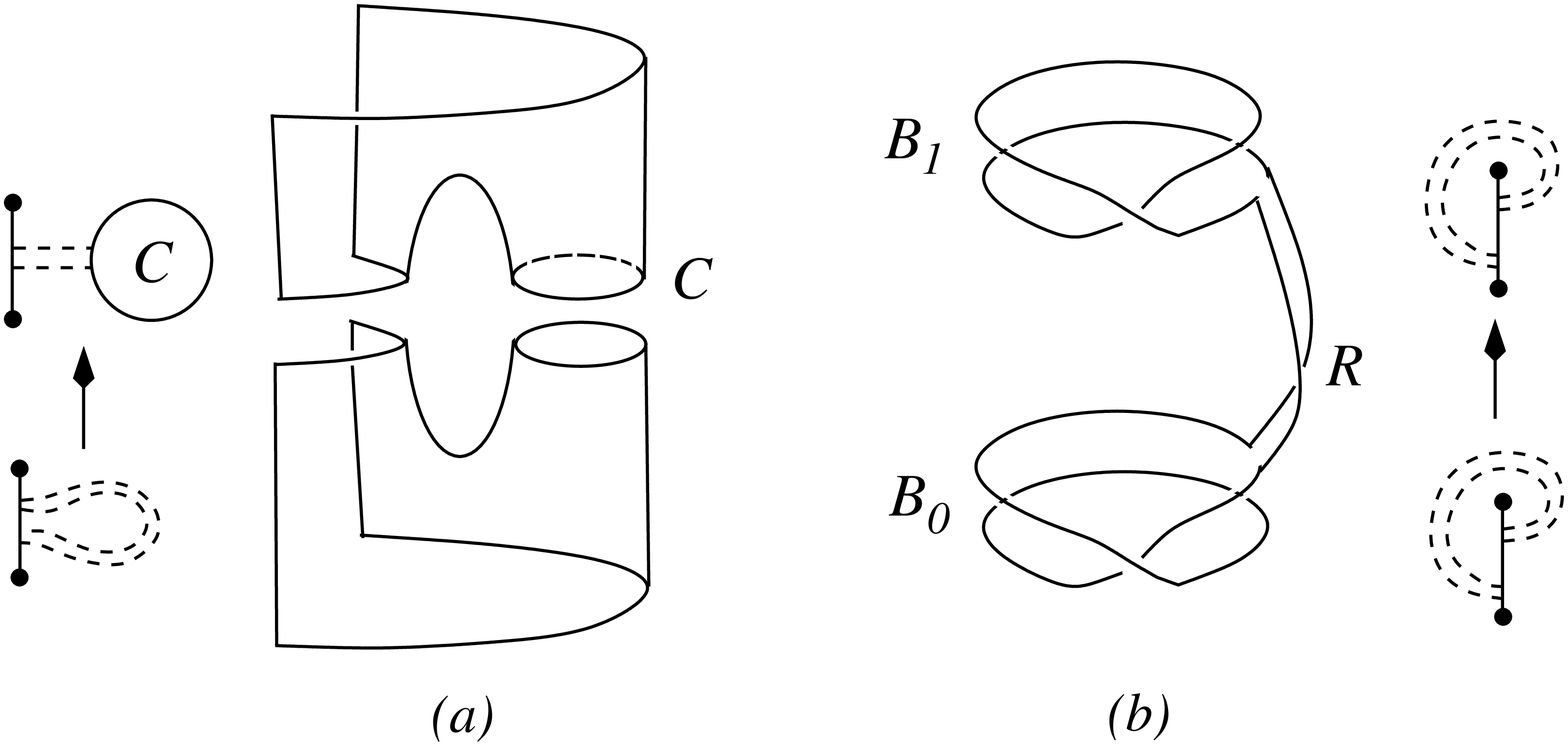}{Construction of $F$ via saddles
in Case (C).}{aa15}
\end{figure}

In this case the saddles, when read from bottom to top and top to
bottom, must join the arcs $S_0\cap F,S_1\cap F$ with themselves,
respectively, both in an orientable fashion or both in a
nonorientable fashion; the possible cases are described (abstractly)
in Fig.\ref{aa15}. In the case of Fig.\ref{aa15}(a), if the level
circle $C$ has slope $(p,q)$ relative to $S_0$, then there is an
essential annulus in $S\times I\setminus\wt F$ with boundary slope
$(p,q)$ in $S_0$. Hence $|p|=1$ or $|q|=1$ since $S\times I$ has no
$F$-spanning annuli, so $K$ is a 2-bridge knot by the argument of
\cite[pp.\ 2161--2162]{matsuda2}.

In the case of Fig.\ref{aa15}(b) let $0<r_1<r_2<1$ be the saddle
levels and, for $i=1,2$, let $B_i$ be the Moebius band $F\cap
S\times[r_i-\varepsilon,r_i+\varepsilon]$ for a sufficiently small
$\varepsilon>0$. Then $F\cap
S\times[r_1+\varepsilon,r_2-\varepsilon]$ is a rectangle $R$, and
$B_1\cup R\cup B_2$ is a once-punctured Klein bottle isotopic to $F$
in $S^3$. Hence $K$ is a knot of the form $K(t_0^*,t_1^*,R')$, where
$t_i$ is the core of the Moebius band $B_i$ in the level $r_i$ and
$R'$ is a monotone subrectangle of $R$.

{\bf Case (D):} \ {\it $F_0'$ is a Moebius band and $F_1'$ is
empty.}

Suppose the first saddle below the 1-level joins the arc component
of $F\cap S_1$ with itself in an orientable fashion; then the first
saddle above the 0-level necessarily joins the arc component of
$F\cap S_0$ with itself in a nonorientable fashion. The situation
here is similar to that of Case (A): the circle $\partial F_0'$
bounds an annulus $A$ in $S\times I$ which can be isotoped away from
$\wt F$ (see Fig.~\ref{aa14}(a), with $C=\partial F_0'$), hence the
slope of $\partial A$ in $S_1$ must be integral and so $K$ is a
2-bridge knot.

Otherwise, the first saddle below the 1-level, say at level
$0<r_1<1$, joins the arc component of $F\cap S_1$ with itself in a
nonorientable fashion, while the first saddle above the 0-level
joins the arc component of $F\cap S_0$ with the circle $\partial
F'_0$. This time the situation is similar to that of Cases (B) and
the second part of (C): for a small $\varepsilon>0$, if $B_1$ is the
Moebius band $F\cap S\times [r_1-\varepsilon,r_1+\varepsilon]$, then
$R=F\cap S\times[0,r_1-\varepsilon]$ is a rectangle and $F_0'\cup
R\cup B_1$ is a once-punctured Klein bottle isotopic to $F$ in
$S^3$, hence $K$ is a knot of the form $K(t_0^*,t_1^*,R')$, where
$R'$ is a monotone subrectangle of $R$ and $t_0,t_1$ can be
described as in Cases (B) and (C), respectively.
\end{proof}

\bibliographystyle{amsplain}

\providecommand{\bysame}{\leavevmode\hbox to3em{\hrulefill}\thinspace}
\providecommand{\MR}{\relax\ifhmode\unskip\space\fi MR }
\providecommand{\MRhref}[2]{%
  \href{http://www.ams.org/mathscinet-getitem?mr=#1}{#2}
}
\providecommand{\href}[2]{#2}

\end{document}